\documentstyle[11pt]{article}
\begin{document}

\title{Surfaces in 3-space possessing nontrivial deformations   \\
which preserve the shape operator}
\author{{\Large Ferapontov E.V. } \\
    Department of Mathematical Sciences \\
    Loughborough University \\
    Loughborough, Leicestershire LE11 3TU \\
    United Kingdom \\
    e-mail: {\tt E.V.Ferapontov@lboro.ac.uk} \\
    and \\
    Centre for Nonlinear Studies \\
    Landau Institute of Theoretical Physics\\
    Academy of Science of Russia, Kosygina 2\\
    117940 Moscow, GSP-1, Russia\\
}
\date{}
\maketitle

\newtheorem{theorem}{Theorem}

\pagestyle{plain}

\maketitle

\begin{abstract}

\bigskip

\noindent The class of surfaces in 3-space possessing nontrivial
deformations 
which preserve principal directions and principal curvatures
(or, equivalently, the shape operator) was investigated by
Finikov and Gambier as far
back as in 1933. We  review some of the known examples and results,
demonstrate the integrability of the corresponding Gauss-Codazzi
equations and draw parallels between this geometrical problem
and the theory of compatible Poisson brackets of hydrodynamic type.
It turns out that  coordinate hypersurfaces of the n-orthogonal systems
arising in the theory of compatible Poisson brackets of hydrodynamic
type must necessarily possess
deformations preserving the shape operator.

\bigskip

\end{abstract}

\newpage
\section{Introduction}
\noindent In 1865 Bonnet initiated the study of the following two
types of deformations of surfaces in the Euclidean 3-space:

\medskip

\noindent 1. {\it Isometric deformations preserving principal directions}
\cite{Bonnet1};

\medskip

\noindent2. {\it Isometric deformations preserving principal curvatures}
\cite{Bonnet2}.

\medskip

\noindent Deformations of the type 1 characterise the so-called moulding
surfaces
(see \cite{Bryant} for further discussion), while deformations of the type 2
lead to  constant mean curvature surfaces and a special class of
Weingarten surfaces, which has been recently investigated
in \cite{Bobenko}. In 1933 Finikov and Gambier \cite{Finikoff},
\cite{Finikoff1}
(see also \cite{Cartan})
introduced deformations of the third type which seem to be a natural
counterpart of deformations
1, 2:

\medskip

\noindent 3. {\it Deformations preserving principal directions and principal
curvatures 
or, equivalently, the shape operator (Weingarten operator)}.

\medskip

\noindent For short, we will call them  $S$-deformations  (notice that
$S$-deformations are
no longer isometric).
Let $M^2\in E^3$ be a surface in the Euclidean space $E^3$
parametrized by coordinates $R^1, R^2$ of curvature lines. Let
\begin{equation}
G_{11}(dR^1)^2+G_{22}(dR^2)^2
\label{1}
\end{equation}
be its third fundamental form (i.e. metric of Gaussian image, which is
automatically of 
constant curvature 1). This condition can be written in the form
\begin{equation}
(\partial_2a+a^2)G^{22}+\frac{a}{2}\partial_2G^{22}+
(\partial_1b+b^2)G^{11}+\frac{b}{2}\partial_1G^{11}+1=0,
\label{K}
\end{equation}
where $\partial_1=\partial / \partial R^1, \
 \partial_2=\partial / \partial R^2, \ G^{11}=1/G_{11}, \  G^{22}=1/G_{22}$
and the coefficients $a$ and $b$ are defined by the formulae
\begin{equation}
\partial_2G^{11}=-2aG^{11}, ~~~  \partial_1G^{22}=-2bG^{22}.
\label{KK}
\end{equation}
Notice that $a$ and $b$ are the Christoffel symbols
of the Levi-Civita connection of the metric (\ref{1}): $a=\Gamma^1_{12}, \
b=\Gamma^2_{21}$. We point out that for given $a$ and $b$ equations
(\ref{K}), (\ref{KK})
are linear in $G^{11}$, $G^{22}$.
Let $k^1, k^2$ be the radii of principal curvature of the surface $M^2$.
In the coordinates of curvature lines Codazzi equations take the form
\begin{equation}
\frac{\partial_2k^1}{k^2-k^1}=a, ~~~~
\frac{\partial_1k^2}{k^1-k^2}=b.
\label{C}
\end{equation}
Let us assume now that the surface $M^2$ possesses nontrivial
$S$-deformation.
Analytically, this means that for given $k^1, k^2$ (and hence given $a, b$)
the
linear system (\ref{K}), (\ref{KK}) is not uniquely solvable for $G^{11},
G^{22}$.
The form of the equations (\ref{K}) --- (\ref{C}) leads to the following
simple observations:

1. $S$-deformations occur in one-parameter (multi-parameter) families.
Indeed,
if $(G^{11}, G^{22})$ and  $(\tilde G^{11}, \tilde G^{22})$ are two
different
solutions of (\ref{K}), (\ref{KK}) (with $a$ and $b$ fixed) then
$(\lambda G^{11}+(1-\lambda) \tilde G^{11}, \
\lambda G^{22}+(1-\lambda) \tilde G^{22})$ is a solution as well. This is a
consequence of the linearity of equations (\ref{K}), (\ref{KK}) in $G^{11},
G^{22}$.

2. The property for a surface $M^2$ to possess nontrivial $S$-deformations
entirely depends on geometry of the orthogonal net on the sphere $S^2$ which
is the Gauss image of its curvature lines.
Indeed, equations (\ref{K}), (\ref{KK}) depend  on $a$ and $b$ only, so that
any
solution of the linear system (\ref{C}) defines the radii of principal
curvature
of an $S$-deformable surface. Surfaces corresponding to different solutions
of (\ref{C}) have one and the same spherical image of curvature lines
 and hence one and the same third fundamental form.
Thus, they $S$-deform simultaneosly. We refer to \cite{Finikoff} for further
comments.

In section 2 we list some of the known examples of $S$-deformable surfaces.
These, in particular, include arbitrary quadrics, cyclids of Dupin and
conformal images
of the surfaces of revolution (as well as all other surfaces with the same
spherical 
image of curvature lines). Apart from several degenerate cases, these
examples exhaust the list
of surfaces which possess 3-parameter families of $S$-deformations
(under certain assumptions the number 3 proves to be the maximal possible).

Following \cite{Fer00}, in section 3 we demontrate the integrability of the
Gauss-Codazzi equations governing
$S$-deformable surfaces by explicitly constructing the Lax pair with a
spectral parameter.

Some results on multi-dimensional $S$-deformable hypersurfaces are discussed
in section 4.
It is pointed out that under certain ``genericity''
assumptions $n$-dimensional quadrics are the only hypersurfaces which
possess $S$-deformations depending on $n+1$ arbitrary constants.

In sections 5---7 we draw parallels between $S$-deformable hypersurfaces and
compatible Poisson brackets of hydrodynamic type. It turns out that
coordinate 
hypersurfaces of the $n$-orthogonal coordinate systems arising in the theory
of bi-Hamiltonian
systems of hydrodynamic type are necessarily $S$-deformable.

\section{Examples}

This section contains a list of examples of surfaces
 possessing nontrivial $S$-deformations.
In our discussion we mainly follow \cite{Finikoff}.

{\bf Example 1. Moulding surfaces of Monge} correspond to the third
fundamental form
\begin{equation}
\frac{\displaystyle (dR^1)^2}{\displaystyle 1+{c}/{\cos ^2 (R^1)}}+
\frac{\displaystyle \sin ^2 (R^1) (dR^2)^2}{\displaystyle \psi (R^2)},
\label{Monge}
\end{equation}
where $c$ is a constant and $\psi (R^2)$ is an arbitrary function. The
metric 
(\ref{Monge}) has constant curvature $1$ for any $c, \ \psi$.
The coordinate curves $R^2=const$ ($R^1=const$) represent the meridians
(parallels) of the sphere $S^2$, respectively.
It  can be readily verified that
 the radii of principal curvature $k^1, k^2$ satisfy the linear system
$$
\partial_2k^1=0, ~~~~ \frac{\partial_1k^2}{k^1-k^2}=\cot (R^1)
$$
which does not explicitly depend on $c, \psi$. Thus, any surface with the
third fundamental form (\ref{Monge}) possesses $S$-deformations depending
on one arbitrary constant and one arbitrary function.
Geometrically, the moulding surfaces of Monge are swept by a planar
curve while the plane rolls along  a curvilinear cylinder.

{\bf Example 2. General moulding surfaces} are characterised by the third
fundamental form
\begin{equation}
(dR^1)^2+
\frac{\displaystyle \sin ^2 (R^1+\varphi(R^2)) (dR^2)^2}{\displaystyle \psi
(R^2)},
\label{moulding}
\end{equation}
where $\psi (R^2)$ and $\varphi(R^2)$ are arbitrary functions.
The coordinate curves $R^1=const$  represent
a one-parameter family of great circles on $S^2$,
the curves $R^2=const$ are their orthogonal trajectories.
It  can be readily verified that
 the radii of principal curvature $k^1, k^2$
satisfy the equations
$$
\partial_2k^1=0, ~~~~ \frac{\partial_1k^2}{k^1-k^2}=\cot (R^1+\varphi(R^2))
$$
which do not explicitly depend on $ \psi$. Thus, any surface with the
third fundamental form (\ref{moulding}) possesses $S$-deformations depending
on one arbitrary function of one variable.
A general moulding surface is swept by a planar
curve while the plane rolls along  a developable surface.

If a surface is not the moulding one
(that is, both $\partial_2 k^1$ and $\partial_1k^2$
are nonzero), the number of parameters on which $S$-deformations  depend
cannot exceed 3. The proof of this statement can be found in \cite{Finikoff}
(see also \cite{Fer90}). Examples presented below actually
exhaust the list of surfaces which possess exactly 3-parameter families of
$S$-deformations.

{\bf Example 3. Quadrics} correspond to the third
fundamental form
\begin{equation}
\frac{R^2-R^1}{4}\left(\frac{\displaystyle (dR^1)^2}
{\displaystyle (R^1)^3+a(R^1)^2+bR^1+c}-
\frac{\displaystyle (dR^2)^2}
{\displaystyle (R^2)^3+a(R^2)^2+bR^2+c}\right)
\label{quadrics}
\end{equation}
where $a, b, c$ are arbitrary constants.
The coordinates $R^1, R^2$ are known as spherical-conical: the corresponding
coordinate curves are
intersections of the sphere with a confocal family of quadratic cones.
Equations for the radii of principal curvature $k^1, k^2$
take the form
$$
\frac{\partial_2k^1}{k^2-k^1}=\frac{1}{2(R^2-R^1)}, ~~~~
\frac{\partial_1k^2}{k^1-k^2}=\frac{1}{2(R^1-R^2)}
$$
which does not explicitly depend on $a, b, c$. Thus, any surface with the
third fundamental form (\ref{quadrics}) possesses $S$-deformations depending
on three arbitrary constants. The case of quadrics corresponds to the choice
$$
k^1=\frac{1}{R^1\sqrt {R^1R^2}}, ~~~ k^2=\frac{1}{R^2\sqrt {R^1R^2}}.
$$

{\bf Example 4. Cyclids of Dupin} correspond to the third
fundamental form
\begin{equation}
\frac{1}{(R^1-R^2)^2}\left(\frac{\displaystyle (dR^1)^2}
{\displaystyle a(R^1)^2+bR^1+c}-
\frac{\displaystyle (dR^2)^2}
{\displaystyle a(R^2)^2+bR^2+c+1}\right)
\label{Dupin}
\end{equation}
where $a, b, c$ are arbitrary constants. The coordinate system $R^1, R^2$
consists of the two
orthogonal pencils of circles on $S^2$.
It can be readily verified that
 the radii of principal curvature $k^1, k^2$
satisfy the linear system
$$
\frac{\partial_2k^1}{k^2-k^1}=\frac{1}{R^1-R^2}, ~~~~
\frac{\partial_1k^2}{k^1-k^2}=\frac{1}{R^2-R^1}
$$
which  does not explicitly depend on $a, b, c$. Thus, any surface with the
third fundamental form (\ref{Dupin}) possesses $S$-deformations  depending
on three arbitrary constants. The case of cyclids of Dupin corresponds to
the choice
$$
k^1=R^2, ~~~ k^2=R^1.
$$

{\bf Example 5. Conformal transforms of  surfaces of revolution} correspond
to the third 
fundamental form
\begin{equation}
\frac{1}{(R^1-R^2)^2}\left(\frac{\displaystyle p^2 (dR^1)^2}
{\displaystyle a(R^1)^2+bR^1+c}-
\frac{\displaystyle (p+p'(R^1-R^2))^2 (dR^2)^2}
{\displaystyle a(R^2)^2+bR^2+c+p^2}\right)
\label{conformal}
\end{equation}
where $a, b, c$ are arbitrary constants and $p$ is an arbitrary function of
$R^2$. 
Here the curves $R^2=const$ represent a one-parameter family of circles on
$S^2$.
It can be readily verified that
 the radii of principal curvature $k^1, k^2$
satisfy the system
$$
\frac{\partial_2k^1}{k^2-k^1}=\frac{1}{R^1-R^2}+\frac{p'}{p}, ~~~~
\frac{\partial_1k^2}{k^1-k^2}=\frac{1}{R^2-R^1}+\frac{p'}{p+p'(R^1-R^2)}
$$
which does not explicitly depend on $a, b, c$. Thus, any surface with the
third fundamental form (\ref{conformal}) possesses $S$-deformations
 depending on three arbitrary constants. The case of conformal images of
surfaces of revolution corresponds to the choice
\begin{equation}
k^1=\frac{1}{p}, ~~~ k^2=\frac{1}{p+p'(R^1-R^2)}
\label{k}
\end{equation}
(see \cite{Fer90}).
The choice $p=R^2$ and the subsequent change of variables
$R^i\to 1/R^i$ returns to the case of cyclids of Dupin discussed above.
Indeed, cyclids of Dupin are known to be  conformal images of cylinders,
cones and tori of revolution. Notice that equations (\ref{k}) possess a more
general class of
solutions
$$
k^1=q, ~~~ k^2=q + \frac{q p'(R^1-R^2)}{p+p'(R^1-R^2)}
$$
where $q$ is another arbitrary function of $R^2$. Solution (\ref{k})
corresponds to
the choice $q=1/p$.

Examples 3-5 provide a complete list of surfaces possessing 3-parameter
families of $S$-deformations.
We refer to \cite{Fer90}
for the proof of this result. Notice that  example 5 is missing
from the classification
proposed in \cite{Finikoff} (although examples 3 and 4 are present).

The case of surfaces possessing  2-parameter families of
$S$-deformations is not undestood at the moment well enough.
We  present here just one example of that type.

{\bf Example 6. Surfaces possessing 2-parameter family of $S$-deformations.}
Consider the metric
\begin{equation}
\frac{R^1-R^2}{(R^1+R^2)^2}\left(\frac{\displaystyle (dR^1)^2}
{\displaystyle a(R^1)^2-R^1+c}-
\frac{\displaystyle (dR^2)^2}
{\displaystyle a(R^2)^2-R^2+c}\right)
\label{2-parameter}
\end{equation}
where $a, c$ are arbitrary constants.
The metric (\ref{2-parameter}) has constant curvature $1$ for any values of
$a, c$.
The coordinate system $R^1, R^2$ is the image of the
ellipsoidal coordinate system on the plane $E^2$ under the stereographic
projection
$E^2\to S^2$ (see \cite{Finikoff}).
Equations for the radii of principal curvature $k^1, k^2$
$$
\frac{\partial_2k^1}{k^2-k^1}=\frac{1}{2(R^2-R^1)}-\frac{1}{R^1+R^2}, ~~~~
\frac{\partial_1k^2}{k^1-k^2}=\frac{1}{2(R^1-R^2)}-\frac{1}{R^1+R^2}
$$
 do not explicitly depend on $a, c$. Thus, any surface with the
third fundamental form (\ref{2-parameter}) possesses $S$-deformations
 depending on two arbitrary constants. The particular case
$$
k^1=\frac{(R^1+R^2)^2}{(R^1)^{5/2}(R^2)^{3/2}}, ~~~
k^2=\frac{(R^1+R^2)^2}{(R^1)^{3/2}(R^2)^{5/2}}
$$
was discussed in detail in \cite{Finikoff}.

{\bf Example 7. Surfaces possessing 1-parameter family of $S$-deformations.}
Consider the metric
\begin{equation}
\frac{2}{\cosh ^2(R^1+R^2)}\left(\frac{\displaystyle (dR^1)^2}
{\displaystyle 1+c}+
\frac{\displaystyle (dR^2)^2}
{1-c}\right)
\label{oneparameter}
\end{equation}
where $c$ is an arbitrary constant. It has constant curvature 1 for any
value of $c$.
The radii of principal curvature $k^1, k^2$
satisfy the equations
$$
\frac{\partial_2k^1}{k^2-k^1}=
\frac{\partial_1k^2}{k^1-k^2}=-\tanh(R^1+R^2)
$$
which do not explicitly depend on $c$. Thus, any surface with the
third fundamental form (\ref{oneparameter}) possesses $S$-deformations
 which  depend
on one arbitrary constant.
The  choice
$$
k^1=-k^2=\cosh^2(R^1+R^2)
$$
corresponds to a minimal surface.
We refer to \cite{Finikoff} for the
geometry of this particular example. The discussion of surfaces possessing
one-parameter
families of $S$-deformations will be continued in the next section.

\section{The Lax pair}

In this section we discuss surfaces which possess 1-parameter
families of $S$-deformations.
Let $M^2\in E^3$ be a surface parametrized by coordinates
$R^1, R^2$ of curvature lines with the third fundamental form (\ref{1}):
$$
G_{11}(dR^1)^2+G_{22}(dR^2)^2.
$$
 Let $k^1, k^2$ be the radii of principal curvature
satisfying the Codazzi equations
\begin{equation}
\frac{\partial_2k^1}{k^2-k^1}=\partial_2\ln \sqrt{G_{11}}, ~~~~
\frac{\partial_1k^2}{k^1-k^2}=\partial_1\ln \sqrt{G_{22}}.
\label{PC}
\end{equation}
Suppose there exists a flat metric
\begin{equation}
g_{11}(dR^1)^2+g_{22}(dR^2)^2
\label{2}
\end{equation}
such that
\begin{equation}
G_{11}=g_{11}/\eta_1, ~~~~ G_{22}=g_{22}/\eta_2,
\label{*}
\end{equation}
where $\eta_1, \eta_2$ are functions of $R^1, R^2$, respectively. One can
readily verify that under these assumptions
the metric
\begin{equation}
\tilde G_{11}(dR^1)^2+\tilde G_{22}(dR^2)^2=
\frac{g_{11}}{\lambda +\eta_1}(dR^1)^2+\frac{g_{22}}{\lambda+\eta_2}(dR^2)^2
\label{3}
\end{equation}
has constant curvature $1$ for any $\lambda$. Since equations (\ref{PC})
are still true if 
we replace $G_{ii}$ by $\tilde G_{ii}$, we arrive at a 1-parameter family of
surfaces $M^2_{\lambda}$
with the third fundamental forms  (\ref{3}) (which depend on $\lambda$) and
the principal curvatures
$k^1, k^2$ (which are independent of $\lambda$). Hence, the shape  operators
of 
 surfaces  $M^2_{\lambda}$ coincide. The problem of the classification of
surfaces which possess
1-parameter families of $S$-deformations  is thus reduced to the
classification
of  metrics (\ref{3}) which have constant
 Gaussian curvature  $1$ for any $\lambda$.
Any such   metric generates an infinite family of $S$-deformable surfaces
whose principal curvatures
$k^1, k^2$ satisfy (\ref{PC}). In terms of the Lame coefficients
$H_1=\sqrt{g_{11}}, \ H_2=\sqrt{g_{22}}$ and
the rotation coefficients $\beta_{12}=\partial_1H_2/H_1, \
\beta_{21}=\partial_2H_1/H_2$ our problem reduces
to the nonlinear system
\begin{equation}
\begin{array}{c}
\partial_1H_2=\beta_{12}H_1, ~~~ \partial_2H_1=\beta_{21}H_2, \\
\ \\
\partial_1\beta_{12}+\partial_2\beta_{21}=0, \\
\ \\
\eta_1 \partial_1\beta_{12}+\eta_2
\partial_2\beta_{21}+\frac{1}{2}\eta_1'\beta_{12}
+\frac{1}{2}\eta_2'\beta_{21}+H_1H_2=0,
\end{array}
\label{4}
\end{equation}
which possesses the Lax pair
$$
\begin{array}{c}
\partial_1\psi=
\left(\begin{array}{ccc}
0 & -\sqrt{\frac{\lambda + \eta_2}{\lambda+\eta_1}}\beta_{21} &
\frac{H_1}{\sqrt{\lambda+\eta_1}}\\
\sqrt{\frac{\lambda + \eta_2}{\lambda+\eta_1}}\beta_{21} & 0 & 0 \\
-\frac{H_1}{\sqrt{\lambda+\eta_1}} & 0 & 0
\end{array}\right)
\psi, \\
\ \\
\partial_2\psi=
\left(\begin{array}{ccc}
0 & \sqrt{\frac{\lambda + \eta_1}{\lambda+\eta_2}}\beta_{12} & 0\\
-\sqrt{\frac{\lambda + \eta_1}{\lambda+\eta_2}}\beta_{12} & 0 &
\frac{H_2}{\sqrt{\lambda+\eta_2}}  \\
0 & -\frac{H_2}{\sqrt{\lambda+\eta_2}} & 0
\end{array}\right)
\psi.
\end{array}
$$
Geometrically, this Lax pair governs infinitesimal displacements of
 the orthonormal frame of the orthogonal coordinate
system on the unit sphere $S^2$, corresponding to the metric (\ref{3}).
In $2\times 2$ matrices it takes the form
$$
\begin{array}{c}
2\sqrt{\lambda+\eta_1}\ \partial_1\psi=
\left(\begin{array}{ccc}
 i\sqrt{{\lambda + \eta_2}}\beta_{21} & {H_1}\\
-{H_1} & - i\sqrt{{\lambda + \eta_2}} \beta_{21}
\end{array}\right) \psi, \\
\ \\
2\sqrt{\lambda+\eta_2}\ \partial_2\psi=
i \left(\begin{array}{ccc}
 -\sqrt{{\lambda + \eta_1}}\beta_{12} & {H_2}\\
{H_2} & \sqrt{{\lambda + \eta_1}} \beta_{12}
\end{array}\right) \psi.
\end{array}
$$

{\bf Example 8. The case of constant $\eta_i$.} Let us assume
 $\eta_1=-1/2, \ \eta_2=1/2$ in which case equations (\ref{4}) take the form
\begin{equation}
\begin{array}{c}
\partial_1H_2=\beta_{12}H_1, ~~~ \partial_2H_1=\beta_{21}H_2, \\
\ \\
\partial_1\beta_{12}=H_1H_2, ~~~~ \partial_2\beta_{21}=-H_1H_2.
\end{array}
\label{40}
\end{equation}
It can be readily verified that system (\ref{40}) possesses 2 integrals
$$
\partial_1(\beta_{12}^2-H_2^2)=0, ~~~~ \partial_2(\beta_{21}^2+H_1^2)=0
$$
implying
$$
\beta_{12}^2-H_2^2=\mu_2(R^2), ~~~~ \beta_{21}^2+H_1^2=\mu_1(R^1),
$$
where the functions $\mu_1, \mu_2$ can be reduced to $\pm 1$ by virtue of
the 
following symmetry of system (\ref{40}):
$$
\begin{array}{c}
\tilde R^1=s_1(R^1), ~~~ \tilde H_1=H_1/s_1', ~~~
\tilde \beta_{21}=\beta_{21}/s_1', \\
\ \\
\tilde R^2=s_2(R^2), ~~~ \tilde H_2=H_2/s_2', ~~~
\tilde \beta_{12}=\beta_{12}/s_2'
\end{array}
$$
(here $s_1(R^1), \ s_2(R^2)$ are arbitrary functions). Let us assume, for
instance, that
$$
\beta_{12}^2-H_2^2=1, ~~~~ \beta_{21}^2+H_1^2=1.
$$
Introducing the parametrisation
$$
H_1=\sin \psi, ~~~ H_2=\sinh \varphi, ~~~
\beta_{12}=\cosh \varphi, ~~~ \beta_{21}=\cos \psi,
$$
we readily rewrite (\ref{40}) in the form
$$
\partial_1 \varphi=\sin \psi, ~~~ \partial_2 \psi = \sinh \varphi,
$$
which implies the following integrable Monge-Ampere equations for $\varphi,
\psi$:
$$
\partial_1\partial_2 \varphi=\sinh \varphi \sqrt{1-(\partial_1 \varphi)^2},
~~~
\partial_1\partial_2 \psi=\sin \psi \sqrt{1+(\partial_2 \psi)^2}.
$$

\section{Multidimensional hypersurfaces  possessing
nontrivial $S$-deformations}

Let $M^n$ be a hypersurface in $E^{n+1}$ parametrised by the
coordinates $R^i$ of curvature lines.
We point out that generic multidimensional hypersurface does not necessarily
possess
curvature line parametrisation. However, (at least for $n=3$)
such parametrisation appears to be the necessary condition for the existence
of  nontrivial $S$-deformations. Let $k^i(R)$ be the radii of
 principal curvature of the hypersurface $M^n$
and $G_{ii}\ (dR^i)^2$ the third
fundamental form 
(which is automatically of constant curvature $1$). These objects satisfy
the Codazzi equations
$$
\frac{\partial_jk^i}{k^j-k^i}=\partial_j \ln \sqrt{G_{ii}}, ~~~ i\ne j.
$$
We will discuss just one important example.

{\bf Example 9. Hyperquadrics} are characterized by the third
fundamental form
\begin{equation}
\sum_{i=1}^n \ \frac{ \prod_{k\ne i}(R^k-R^i)}{ 4 P(R^i)} \ (dR^i)^2,
\label{hyperquadric}
\end{equation}
where $P(R)=\prod_{s=1}^{n+1}(R-a_s)$ is an arbitrary polynomial of the
order $n+1$.
The corresponding radii of principal curvature $k^i$  satisfy the linear
system
$$
\frac{\partial_jk^i}{k^j-k^i}=\frac{1}{2(R^j-R^i)}, ~~~ i\ne j,
$$
which does not explicitly depend on $a_s$. Thus, any hypersurface $M^n$ with
the third fundamental form
(\ref{hyperquadric}) possesses $(n+1)$-parameter family of $S$-deformations.
The case of a hyperquadric corresponds to the choice
$$
k^i=\frac{1}{ R^i\sqrt{\prod_1^nR^k}}.
$$

It can be demonstrated that quadrics are the only $n$-dimensional
hypersurfaces
which possess $(n+1)$-parameter families of $S$-deformations
and satisfy the 
``genericity'' assumption $\partial_j k^i\ne 0$ for any $i\ne j$
(we recall that 
for $n=2$ the condition
$\partial_1k^2 \ne 0, \ \partial_2 k^1\ne 0$ forbids moulding surfaces).
Thus, unlike the case $n=2$, the class of multidimensional hypersurfaces
possessing
the ``maximal'' number of $S$-deformations
appears to be very restricted. It is likely that there exist
examples of hypersurfaces possessing $k$-parameter families of
$S$-deformations 
for any intermediate $k=1, 2, ..., n+1$. Numerous examples of
that type are provided by coordinate hypersurfaces of the n-orthogonal
coordinate systems 
arising in the theory of multi-Hamiltonian equations of hydrodynamic type.

\section{Compatible Poisson brackets of hydrodynamic type.
Criterion of the compatibility.}

In 1983 Dubrovin and Novikov \cite{DN83} introduced the Poisson brackets of
hydrodynamic type
\begin{equation}
\{F, G\}=\int \frac{\delta F}{\delta u^i} A^{ij} \frac{\delta G}{\delta u^j}
\ dx
\label{Poisson}
\end{equation}
defined by the Hamiltonian operators $A^{ij}$ of the form
\begin{equation}
A^{ij}=g^{ij}\frac{d}{dx} +b^{ij}_ku^k_x, ~~~~ b^{ij}_k=-g^{is}\Gamma
^j_{sk}.
\label{H1}
\end{equation}
They proved that in the nondegenerate case $(det ~g^{ij}\ne 0)$ the bracket
(\ref{Poisson}), (\ref{H1}) is
skew-symmetric and satisfies the Jacobi identities if and only if the metric
$g^{ij}$ (with upper indices) is flat, and $\Gamma^j_{sk}$
are the Christoffel symbols of the corresponding Levi-Civita connection.

Let us consider another Poisson bracket of hydrodynamic type defined on the
same phase
space by the Hamiltonian operator
\begin{equation}
\tilde A^{ij}=\tilde g^{ij}\frac{d}{dx}+\tilde b^{ij}_ku^k_x, ~~~~
\tilde b^{ij}_k=-\tilde g^{is}\tilde \Gamma ^j_{sk},
\label{H2}
\end{equation}
corresponding to a flat metric $\tilde g^{ij}$. Two Poisson brackets
(Hamiltonian operators) are called compatible if their linear combinations
$\tilde A^{ij}+\lambda A^{ij}$ are Hamiltonian as well.
This requirement implies, in particular, that the metric
$\tilde g^{ij}+\lambda g^{ij}$ must be flat for any  $\lambda$
(plus certain additional restrictions). The necessary and sufficient
conditions of the compatibility
were first formulated by Dubrovin \cite{D93}, \cite{D94} (see
\cite{Mokhov99}, \cite{Mokhov00} for further discussion).
Below we reformulate these conditions in terms
of the operator 
$
r^i_j=\tilde g^{is} g_{sj}
$ 
(Theorem 1) which, in particular, imply the vanishing of the Nijenhuis
tensor of the operator $r^i_j$:
$$
N^i_{jk}=r^s_j\partial_sr^i_k  -  r^s_k\partial_sr^i_j  -
r^i_s(\partial_j r^s_k - \partial_kr^s_j)=0
$$ 
(see \cite{Fer95}, \cite{Mokhov00}).

Examples of compatible Hamiltonian pairs naturally arise in the theory of
Hamiltonian systems of hydrodynamic type --- see e.g.
\cite{Arik}, \cite{Gumral90}, \cite{Gumral94}, \cite{Nutku}, \cite{Olver},
 \cite{Pavlov94.1}.
Dubrovin developed a deep theory for the particular class of compatible
Poisson
 brackets  arising in the framework of the associativity equations
\cite{D93}, \cite{D94}. Compatible Poisson brackets of
hydrodynamic type can also be obtained as a result of  Whitham
averaging (dispersionless limit) from the local
compatible Poisson brackets of integrable systems
\cite{DN83}, \cite{DN84}, \cite{DN89}, \cite{Tsarev90}, \cite{Pavlov94},
\cite{FerPav},  \cite{Strachan}.
Some further examples and partial classification results  can be found in
\cite{Fer90}, \cite{FerPav}, \cite{Nutku}, \cite{Pavlov93},
\cite{Mokhov99}, \cite{Fordy}.

If the spectrum of $r^i_j$ is simple, the vanishing of the Nijenhuis tensor
 implies the existence of a coordinate system where both metrics $g^{ij}$
and 
$\tilde g^{ij}$ become diagonal. In these  coordinates the compatibility
conditions 
take the form of an integrable reduction of the Lam\'e equations. We present
the 
corresponding Lax pairs in section 6. Another approach to the
integrability of this system was recently proposed by
Mokhov \cite{Mokhov001} by an appropriate modification of Zakharov's scheme
\cite{Zakharov}. 

Our main observation is the relationship between
compatible Poisson brackets of hydrodynamic type and
hypersurfaces $M^{n-1}\in E^n$  possessing nontrivial $S$-deformations.
In section 7 we demonstrate that the n-orthogonal coordinate systems in
$E^n$ 
corresponding to  flat metrics $\tilde g^{ij}+\lambda g^{ij}$
(rewritten in the diagonal coordinates)
deform with respect to $\lambda$ in such a way that the shape operators of
coordinate hypersurfaces are preserved up to  constant scaling factors.

Notice that the operator $r^i_j=\tilde g^{is}g_{sj}$ is automatically
symmetric:
\begin{equation}
r^i_sg^{sj}=r^j_sg^{si},
\label{symm}
\end{equation}
so that $\tilde g^{ij}=r^i_sg^{sj}=r^j_sg^{si}=r^{ij}$. In what follows
we use the first metric $g^{ij}$ for raising and lowering the indices.

\begin{theorem} {\rm \cite{Fer00}}

 Hamiltonian operators (\ref{H1}), (\ref{H2}) are compatible if and
only if the following conditions are satisfied:

\medskip

1. The Nijenhuis tensor of $r^i_j$ vanishes:
\begin{equation}
N^i_{jk}=r^s_j\partial_sr^i_k  -  r^s_k\partial_sr^i_j  -
r^i_s(\partial_j r^s_k - \partial_kr^s_j)=0.
\label{Nijenhuis}
\end{equation}

\medskip

2. The metric coefficients $\tilde g^{ij}= r^{ij}$ satisfy the equations
\begin{equation}
\nabla^i\nabla^j r^{kl}+\nabla^k\nabla^l r^{ij}=
\nabla^i\nabla^k r^{jl}+\nabla^j\nabla^l r^{ik}.
\label{nabla}
\end{equation}
Here $\nabla^i=g^{is}\nabla_s$ is the covariant differentiation
corresponding to the metric
$g^{ij}$. The vanishing of the Nijenhuis tensor implies the following
expression for the
coefficients $\tilde b^{ij}_k$ in terms of $r^i_j$:
\begin{equation}
2 \tilde b^{ij}_k= \nabla^ir^j_k-\nabla^jr^i_k+\nabla_k
r^{ij}+2b^{sj}_kr^i_s
\label{F}
\end{equation}

\end{theorem}

In a somewhat different form the necessary and sufficient conditions of the
compatibility were formulated in \cite{D93}, \cite{D94}, \cite{Mokhov99},
\cite{Mokhov00}.

{\bf Remark.} The criterion of the compatibility of  Hamiltonian operators
of
hydrodynamic type resembles that of  finite-dimensional Poisson bivectors:
two skew-symmetric Poisson bivectors $\omega^{ij}$ and $\tilde \omega^{ij}$
are compatible if and only if the Nijenhuis tensor of the corresponding
recursion operator $r^i_j=\tilde \omega^{is}\omega_{sj}$ vanishes. We
emphasize that
in our situation operator $r^i_j$ does not coincide with the recursion
operator.

\bigskip

{\bf Remark.} If the spectrum of $r^i_j$ is simple the condition
(\ref{nabla}) is redundant:
it is automatically satisfied by virtue of (\ref{Nijenhuis}) and the
flatness of both metrics
$g$ and $\tilde g$. This was the motivation for me to drop condition
(\ref{nabla})
in the compatibility criterion formulated in \cite{Fer95}. However, in this
 general form the criterion
proved to be incorrect:  recently it was pointed out by Mokhov
\cite{Mokhov00} 
that in the case when the spectrum of $r^i_j$ is not simple the vanishing of
the Nijenhuis tensor
is  no longer sufficient for the compatibility.

\section{Compatibility conditions in the diagonal form: the Lax pairs}

If the spectrum of $r^i_j$ is simple, the vanishing of the
Nijenhuis tensor implies the existence of the coordinates $R^1,\ldots, R^n$
in which the  objects $r^i_j, ~ g^{ij}, ~ \tilde g^{ij}$ become diagonal.
Moreover, the $i$-th eigenvalue of $r^i_j$ depends only on the coordinate
$R^i$, so that $r^i_j = diag(\eta_i), ~ g^{ij}= diag(g^{ii}), ~ \tilde
g^{ij} =
diag(g^{ii}\eta_i)$ where $\eta_i$ is a function of $R^i$. This is a
generalization of 
the analogous observation by Dubrovin
\cite{D93} in the particular case of compatible Poisson brackets originating
from the 
theory of the associativity equations. Introducing the Lame coefficients
$H_i$ and the rotation coefficients
$\beta_{ij}$ by the formulae
\begin{equation}
H_i=\sqrt{g_{ii}}=1/\sqrt{g^{ii}}, ~~~~ \partial_i H_j=\beta_{ij} H_i,
\label{H}
\end{equation}
we can rewrite the zero curvature conditions for the metric
$g$ in the form
\begin{equation}
\partial_k\beta_{ij}=\beta_{ik}\beta_{kj},
\label{F1}
\end{equation}
\begin{equation}
\partial_i\beta_{ij}+\partial_j\beta_{ji}+\sum_{k\ne i,
j}\beta_{ki}\beta_{kj}=0.
\label{F2}
\end{equation}
The zero curvature condition for the metric $\tilde g$
imposes the additional constraint
\begin{equation}
\eta_i\partial_i\beta_{ij}+\eta_j\partial_j\beta_{ji}+\frac{1}{2}\eta_i{'}\b
eta_{ij}
+\frac{1}{2}\eta_j{'}\beta_{ji}+\sum_{k\ne i,
j}\eta_k\beta_{ki}\beta_{kj}=0,
\label{F3}
\end{equation}
resulting from (\ref{F2}) after the  substitution of
 the rotation coefficients $\tilde
\beta_{ij}=\beta_{ij}\sqrt{\eta_i/\eta_j}$
of the metric
$\tilde g$. As can be readily seen, equations (\ref{F2}) and (\ref{F3})
already imply the compatibility,
so that in the diagonalisable case condition (\ref{nabla}) of Theorem 1 is
indeed superfluous.
Solving equations (\ref{F2}), (\ref{F3}) for $\partial_i\beta_{ij}$,
 we can rewrite (\ref{F1})--(\ref{F3}) in the form
\begin{equation}
\begin{array}{c}
\partial_k\beta_{ij}=\beta_{ik}\beta_{kj}, \\
\ \\
\partial_i\beta_{ij}=\frac{1}{2} \frac{\displaystyle
\eta_i{'}}{\displaystyle \eta_j-\eta_i}\beta_{ij}
+\frac{1}{2} \frac{\displaystyle \eta_j{'}}{\displaystyle
\eta_j-\eta_i}\beta_{ji}+\sum_{k\ne i, j}
\frac{\displaystyle \eta_k-\eta_j}{\displaystyle
\eta_j-\eta_i}\beta_{ki}\beta_{kj}.
\end{array}
\label{S}
\end{equation}
It can be  verified by a straightforward calculation that system (\ref{S})
is compatible for any choice of the functions
$\eta_i(R^i)$,
and its general solution depends on $n(n-1)$ arbitrary functions of one
variable
(indeed, one can arbitrarily prescribe the value of $\beta_{ij}$ on the j-th
coordinate line). 
Under the additional "Egorov" assumption $\beta_{ij}=\beta_{ji}$, system
(\ref{S})
reduces to the one studied by Dubrovin in \cite{D94}.
For $n\geq 3$ system (\ref{S}) is essentially nonlinear. Its integrability
follows from 
the Lax pair \cite{Fer00}
\begin{equation}
\partial_j\psi_i=\beta_{ij}\psi_j, ~~~~
\partial_i\psi_i=-\frac{\eta_i{'}}{2(\lambda+\eta_i)}\psi_i-\sum_{k\ne
i}\frac{\lambda+\eta_k}
{\lambda+\eta_i}\beta_{ki}\psi_k
\label{L1}
\end{equation}
with  a spectral parameter $\lambda$ (another demonstration of the
integrability of system (\ref{L1})
has been proposed recently in \cite{Mokhov001} by an appropriate
modification of 
Zakharov's approach \cite{Zakharov}).

{\bf Remark.} In fact, the Lax pair (\ref{L1}) is gauge-equivalent to the
equations
$$
({\tilde g}^{ik}+\lambda g^{ik}) \partial _k \partial_j \psi+
({\tilde b}^{ik}_j+\lambda b^{ik}_j) \partial _k \psi=0
$$
for the Casimirs $\int \psi dx$ of the Hamiltonian operator
$\tilde A^{ij}+\lambda A^{ij}$.

After the gauge transformation $\psi_i= \varphi_i/\sqrt{\lambda+\eta_i}$
the Lax pair (\ref{L1}) assumes the manifestly skew-symmetric form
\begin{equation}
\partial _j\varphi _i=\sqrt{\frac{\lambda+\eta _i}{\lambda+\eta _j}}
\beta _{ij}\varphi _j, ~~~~
\partial _i\varphi _i=-\sum_{k\ne i}\sqrt{\frac{\lambda+\eta _k}
{\lambda+\eta _i}}\beta _{ki}\varphi _k.
\label{L2}
\end{equation}
Thus, we can introduce an orthonormal frame $\vec \varphi_1, ..., \vec
\varphi_n$ 
in the Euclidean space $E^n$ satisfying the equations
\begin{equation}
\partial _j\vec \varphi _i=
\sqrt{\frac{\lambda+\eta _i}{\lambda+\eta _j}}\beta_{ij}\vec \varphi _j,
~~~~
\partial _i\vec \varphi _i=-\sum_{k\ne i}\sqrt{\frac{\lambda+\eta _k}
{\lambda+\eta _i}}\beta _{ki}\vec \varphi _k,
~~~~
(\vec \varphi_i, \vec \varphi_j)=\delta_{ij}.
\label{frame}
\end{equation}
Let us introduce a vector $\vec r$ such that
$$
\partial_i \vec r=\frac{H_i}{\sqrt{\lambda+\eta_i}} \vec \varphi_i
$$
(the compatibility of these equations can be readily verified).
In view of the formula
$$
(\partial_i \vec r, \partial_j \vec
r)=\frac{H_i^2}{\lambda+\eta_i}\delta_{ij}
$$
the radius-vector $\vec r$ is descriptive of the n-orthogonal coordinate
system in $E^n$
corresponding to the flat metric
$$
\sum_i \frac{H_i^2}{\lambda+\eta_i}(dR^i)^2.
$$
Geometrically, $\vec \varphi_i$ are the unit vectors along the coordinate
lines of 
this n-orthogonal system.

Let us discuss in some more detail the
 case $\eta_i=const=c_i$, in which system (\ref{S}) takes the form
\begin{equation}
\begin{array}{c}
\partial_k\beta_{ij}=\beta_{ik}\beta_{kj}, \\
\ \\
\partial_i\beta_{ij}=\sum_{k\ne i, j}
\frac{\displaystyle c_k-c_j}{\displaystyle c_j-c_i}\beta_{ki}\beta_{kj}.
\end{array}
\label{S1}
\end{equation}
One can readily verify that the quantity
$$
P_i=\sum_{k\ne i}(c_k-c_i)\beta_{ki}^2
$$
is an integral of system (\ref{S1}), namely,
$\partial_jP_i=0$ for  any $i\ne j$,
so that $P_i$ is a function of $R^i$. Utilising the obvious symmetry
$
R^i\to  s_i(R^i), ~~ \beta_{ki}\to \beta_{ki}/s_i'(R^i)
$
of system (\ref{S1}), we can reduce $P_i$ to $\pm 1$ (if nonzero). Let us
consider the simplest nontrivial case
$n=3, \ P_1=P_2=1, \ P_3=-1$:
$$
\begin{array}{c}
P_1=(c_2-c_1)\beta_{21}^2+(c_3-c_1)\beta_{31}^2=1, \\
\ \\
P_2=(c_1-c_2)\beta_{12}^2+(c_3-c_2)\beta_{32}^2=1, \\
\ \\
P_3=(c_1-c_3)\beta_{13}^2+(c_2-c_3)\beta_{23}^2=-1.
\end{array}
$$
Assuming $c_3>c_2>c_1$ and introducing the parametrisation
$$
\begin{array}{c}
\beta_{21}={\sin p}/{\sqrt {c_2-c_1}}, ~~~~ \beta_{31}={\cos p}/{\sqrt
{c_3-c_1}}, \\
\ \\
\beta_{12}={\sinh q}/{\sqrt {c_2-c_1}}, ~~~~ \beta_{32}={\cosh q}/{\sqrt
{c_3-c_2}}, \\
\ \\
\beta_{13}={\sin r}/{\sqrt {c_3-c_1}}, ~~~~ \beta_{23}={\cos r}/{\sqrt
{c_3-c_2}},
\end{array}
$$
we readily rewrite (\ref{S1}) in the form
$$
\begin{array}{c}
\partial_1q=\mu_1 \cos p, ~~~~ \partial_1r=-\mu_1 \sin p, \\
\ \\
\partial_2p=-\mu_2 \cosh q, ~~~~ \partial_2r=\mu_2 \sinh q, \\
\ \\
\partial_3p=\mu_3 \cos r, ~~~~ \partial_3q=\mu_3 \sin r,
\end{array}
$$
where
$$
\mu_1=\sqrt{\frac{c_3-c_2}{(c_2-c_1)(c_3-c_1)}}, ~~~
\mu_2=\sqrt{\frac{c_3-c_1}{(c_2-c_1)(c_3-c_2)}}, ~~~
\mu_3=\sqrt{\frac{c_2-c_1}{(c_3-c_1)(c_3-c_2)}}.
$$
After rescaling, this system simplifies to
\begin{equation}
\begin{array}{c}
\partial_1q= \cos p, ~~~~ \partial_1r=- \sin p, \\
\ \\
\partial_2p=- \cosh q, ~~~~ \partial_2r= \sinh q, \\
\ \\
\partial_3p= \cos r, ~~~~ \partial_3q= \sin r.
\end{array}
\label{S2}
\end{equation}
Expressing $p$ and $r$ in the form
$
p=\arccos {\partial_1q}, ~ r=\arcsin {\partial_3q},
$
we can rewrite (\ref{S2}) as a triple of pairwise commuting Monge-Ampere
equations
$$
\begin{array}{c}
\partial_1 \partial_2q=\cosh q\sqrt{1-\partial_1q^2}, \\
\ \\
\partial_1 \partial_3q=-\sqrt{1-\partial_1q^2} \sqrt{1-\partial_3q^2}, \\
\ \\
\partial_2 \partial_3q=\sinh q\sqrt{1-\partial_3q^2}.
\end{array}
$$
Similar triples of Monge-Ampere equations  were obtained in \cite{Fer97} in
the classification of quadruples
of $3\times 3$ hydrodynamic type systems which are closed under the Laplace
transformations.
However,  there is no understanding of this  coincidence at the moment.
Notice that coordinate surfaces of these 3-orthogonal coordinate systems are
of the type discussed in example 8.

\section{Deformations of n-orthogonal coordinate systems inducing rescalings
of shape operators of the coordinate hypersurfaces}

We have demonstrated in sect. 6 that the radius-vector $\vec r(R^1,...,
R^n)$
of the n-orthogonal coordinate system in $E^n$ corresponding to the flat
diagonal metric
$ \sum_i \frac{H_i^2}{\lambda+\eta_i}(dR^i)^2$ satisfies the equations
$$
\partial_i \vec r=\frac{H_i}{\sqrt{\lambda+\eta_i}} \vec \varphi_i,
$$
where the infinitesimal displacements of the orthonormal frame $\vec
\varphi_i$
are governed by 
$$
\partial_j\vec \varphi_i=
\sqrt{\frac{\displaystyle \lambda+\eta_i}{\displaystyle
\lambda+\eta_j}}\beta_{ij}\vec \varphi_j, ~~~~
\partial_i\vec \varphi_i=-\sum_{k\ne i}\sqrt{\frac{\displaystyle
\lambda+\eta_k}
{\displaystyle \lambda+\eta_i}}\beta_{ki}\vec\varphi_k.
$$
Since our formulae depend on the spectral parameter, we may speak of the
"deformation" of the n-orthogonal
coordinate system with respect to $\lambda$. To investigate this deformation
in some more detail,
we fix a coordinate hypersurface $M^{n-1}\subset E^n$ (say, $R^n=const$).
Its radius-vector $\vec r$
and the unit normal $\vec \varphi_n$ satisfy the Weingarten equations
$$
\partial_i\vec \varphi_n=\frac{\beta_{ni}}{H_i}\sqrt{\lambda+\eta_n}\
\partial_i \vec r, 
~~~~ i=1,..., n-1,
$$
implying that
$$
k^i=\frac{\beta_{ni}}{H_i}\sqrt{\lambda+\eta_n}
$$
are principal curvatures of the hypersurface  $M^{n-1}$. Since $\eta_n$ is
constant on $M^{n-1}$, our deformation
 preserves the shape operator of $M^{n-1}$ up to a
constant scaling factor $\sqrt{\lambda+\eta_n}$ (we point out that the
curvature line parametrisation
$R^1,..., R^{n-1}$ is preserved by a construction). Thus, compatible Poisson
brackets of hydrodynamic type give rise to deformations of
n-orthogonal systems in $E^n$ which, up to scaling factors,
preserve shape operators of
the coordinate hypersurfaces. If we follow the evolution of a particular
coordinate hypersurface $M^{n-1}$,
this scailing factor can be eliminated by a homothetic transformation of the
ambient space $E^n$,
so that we arrive at the nontrivial deformation which preserves the
shape operator. However, this scaling factor cannot be eliminated for all
coordinate hypersurfaces
simultaneously. 

\section{Concluding remarks}

It seems to be an interesting and nontrivial problem to classify
hypersurfaces 
$M^n\in E^{n+1}$ possessing $S$-deformations which depend
on $k$ essential parameters where
$k$ varies from $2$ to $n$ (the limiting cases $k=n+1$ and $k=1$
reduce to hyperquadrics and generic $S$-deformable surfaces, respectively).
Even the simplest nontrivial case $n=2$, $k=2$ is not understood at the
moment.

\end{document}